\documentclass{amsart}
\usepackage{latexsym}
\usepackage{amssymb}
\usepackage{amsmath}
\usepackage{amsfonts}
\usepackage{amsthm}
\newtheorem{theo}{Theorem}[section]
\newtheorem{pro}[theo]{Proposition}

\newtheorem{cor}[theo]{Corollary}
\newtheorem{defi}[theo]{Definition}
\newtheorem{rem}[theo]{Remark}

\newtheorem*{unnumberedtheorem}{Theorem}

\newcommand{\rn}{\mathbb{R}^{n+1}}

\title{A formula relating entropy monotonicity to Harnack inequalities}
\author{Klaus Ecker}
\date{February 26, 2007}

\begin{document}

\maketitle

\section{Introduction}\label{intro}

In \cite{P}, Perelman considered the functional

\begin{equation*}
\mathcal{W}(g, f, \tau)=\int_X\left(\tau(|\nabla f|^2+R)+f
-(n+1)\right)u\,dV
\end{equation*}
for $\tau >0$ and smooth functions $f$ on a closed $(n+1)$ -
dimensional Riemannian manifold $(X, g)$ where
\begin{equation*}
u=\frac{e^{-f}}{(4\pi\tau)^{\frac{n+1}{2}}}
\end{equation*}
and defined an associated entropy by
\begin{equation*}
\mu (g, \tau)=\inf\left\{\mathcal{W}(g, f, \tau),\,\int_X u
\,dV=1\right\}.
\end{equation*}
His ingenious realization was that when $\tau(t) >0$ satisfies
$\frac{\partial\tau}{\partial t}=-1$, $(X, g(t))$ evolves by the
Ricci flow
\begin{equation*}
\frac{\partial}{\partial t}g_{ij}=-2R_{ij}
\end{equation*}
and $f$ satisfies the equation
\begin{equation*}
\frac{\partial f}{\partial t} +\Delta f +R= |\nabla
f|^2+\frac{n+1}{2\tau}
\end{equation*}
which preserves the condition
\begin{equation*}
\int_X u \,dV=1
\end{equation*}
then
\begin{equation*}
\frac{d}{dt}\mathcal{W}(g(t), f(t), \tau(t))
=2\tau \int_X\left|R_{ij}+\nabla_i\nabla_j f-\frac{g_{ij}}{2\tau}\,\right|^2u\,dV.\nonumber\\
\end{equation*}
This implies in particular that
\begin{equation*}
\frac{d}{dt}\mu (g(t), \tau (t))\ge 0
\end{equation*}
with equality exactly for homothetically shrinking solutions of Ricci flow.

An important consequence of this entropy formula is a lower volume ratio bound for solutions of Ricci flow on a closed manifold for a finite time interval $[0,
T)$ asserting the existence of a constant $\kappa > 0$, only depending on $n, T$ and $g(0)$, such that the inequality
\begin{equation*}
\frac{V_t(B^t_r(x_0))}{r^{n+1}}\ge\kappa
\end{equation*}
holds for all $t\in [0, T)$ and $r\in [0, \sqrt{T})$ for balls $B^t_r(x_0)$ (with respect to $g(t)$) in which the inequality $r^2|Rm|\le 1$ for the Riemann
tensor of $g(t)$ holds.

This lower volume ratio bound rules out certain collapsed metrics as rescaling limits near singularities of Ricci flow such as products of Euclidean spaces with
the so-called cigar soliton solution of Ricci flow given by $X=\mathbb{R}^2$ with the metric
\begin{equation*}
ds^2=\frac{dx^2+dy^2}{1+x^2+y^2}.
\end{equation*}

In this paper, we aim at adapting Perelman's entropy formula to
the situation where a family of bounded open regions
$(\Omega_t)_{t\in [0, T)}$ in $\rn$ with smooth boundary
hypersurfaces $M_t=\partial\Omega_t$ is evolving with smooth
normal speed
\begin{equation*}
\beta_{M_t}=-\frac{\partial x}{\partial t}\cdot\nu.
\end{equation*}
Here $x$ denotes the embedding map of $M_t$ and $\nu$ is the
normal pointing out of $\Omega_t$.

For open subsets $\Omega\subset\rn$ , smooth functions
$f:\bar\Omega\to\mathbb{R}$ and
$\beta:\partial\Omega\to\mathbb{R}$ and $\tau >0$ we consider the
quantity
\begin{equation*}
\mathcal{W}_\beta(\Omega, f, \tau)=\int_\Omega\left(\tau |\nabla
f|^2+f-(n+1)\right)\,u\,dx+2\tau\int_{\partial\Omega}\beta u \,dS
\end{equation*}
with
\begin{equation*}
u=\frac{e^{-f}}{(4\pi\tau)^{\frac{n+1}{2}}}
\end{equation*}
and the associated entropy
\begin{equation*}
\mu_\beta(\Omega, \tau)= \inf \left\{ \mathcal{W}_\beta (\Omega,
f, \tau)\, ,\int_\Omega u\,dx =1 \right\}.
\end{equation*}

We then derive a formula which states that if $(\Omega_t)$ evolves
as above, $\tau(t)
>0$ satisfies $\frac{\partial\tau}{\partial t}=-1$, f satisfies
the evolution equation
\begin{equation*}
\frac{\partial f}{\partial t} +\Delta f = |\nabla
f|^2+\frac{n+1}{2\tau}
\end{equation*}
in $\Omega_t$ with Neumann boundary condition
\begin{equation*}
\nabla f\cdot\nu =\beta
\end{equation*}
on $M_t=\partial{\Omega_t}$ and we introduce a family of
diffeomorphisms  $\varphi_t:\bar\Omega\to\bar\Omega_t$ with
$x=\varphi_t(q),\, q\in\bar\Omega$ obeying
\begin{equation*}
\frac{\partial x}{\partial t}=-\nabla f(x,t)
\end{equation*}
then
\begin{equation*}
\frac{d}{dt}\mathcal{W}_\beta (\Omega_t, f(t), \tau(t)) =2\tau
\int_{\Omega_t}\left|\nabla_i\nabla_j
f-\frac{\delta_{ij}}{2\tau}\,\right|^2u\,dx
 -\int_{M_t}\nabla W\cdot\nu\,dS\nonumber\\
\end{equation*}
where $W=\tau(2\Delta f-|\nabla f|^2)+f-(n+1)$.

For evolving bounded regions $\Omega_t$ inside a fixed Riemannian
manifold $(X, g)$ or inside a Ricci flow solutions one can derive
analoguous versions of this formula.

The main observation in this paper is that this can be converted
to
\begin{align}
\frac{d}{dt}\mathcal{W}_\beta &(\Omega_t, f(t), \tau(t))
=2\tau \int_{\Omega_t}\left| \nabla_i\nabla_j f-\frac{\delta_{ij}}{2\tau}\,\right|^2u\,dx\nonumber\\
& +2\tau\int_{M_t}\left(\frac{\partial\beta}{\partial
t}-2\,\nabla^{M}\beta \cdot \nabla^{M}f
+A(\nabla^{M}f,\nabla^{M}f)-\frac{\beta}{2\tau}\right)u\,dS
\nonumber\end{align} where $A$ denotes the second fundamental form
of $M_t$.

For functions $\beta$ for which the hypersurface integral is
nonnegative the inequality
\begin{equation*}
\frac{d}{dt}\mathcal{W}_\beta (\Omega_t, f(t), \tau(t))
\ge 2\tau \int_{\Omega_t}\left| \nabla_i\nabla_j f-\frac{\delta_{ij}}{2\tau}\,\right|^2u\,dx\nonumber\\
\end{equation*}
results. When $\beta =0$, that is for a fixed bounded region
$\Omega$ with smooth convex boundary inside a fixed manifold of
non-negative Ricci curvature, Lei Ni \cite{N} has previously
obtained this inequality.

It implies, as in Perelman's situation,
\begin{equation*}
\frac{d}{dt}\mu_\beta (\Omega_t, \tau (t))\ge 0
\end{equation*}
and also the following localised lower volume ratio bound:

There is a constant $\kappa >0$ depending only on $n, \Omega_0,
T,\,\sup_{M_0}|\beta|$ and $c_1$ such that
\begin{equation*}
\frac{V(\Omega_t\cap B_r(x_0))}{r^{n+1}}\ge\kappa
\end{equation*}
holds for all $t\in [0, T)$ and $r\in (0, \sqrt {T}]$ in balls $B_r(x_0)\subset\rn$ satisfying the conditions $V(\Omega_t\cap B_{r/2}(x_0))>0$ and
\begin{equation*}
\frac{V(\Omega_t\cap B_r(x_0))+r^2\int_{M_t\cap
B_r(x_0)}|\beta|\,dS}{V(\Omega_t\cap B_{r/2}(x_0))}\le c_1.
\end{equation*}
Since this statement is scaling invariant for suitably homogeneous
$\beta$ it is also valid on any
smooth limit of suitably rescaled solutions of the flow consisting
of smooth, compact embedded hypersurfaces, but now for all radii
$r>0$ as long as the other conditions still hold for the balls
$B_r(x_0)$ we consider.

In the important case of mean curvature flow, that is where
$\beta_{M_t}$ is the mean curvature $H_{M_t}$ of the hypersurfaces
$M_t$, the right hand side of the formula vanishes on
homothetically shrinking solutions and for $f=|x|^2/4\tau$. This
leads us to the following conjecture:

\bigskip

{\bf Conjecture.} {\sl In the case of mean curvature flow in $\rn$
for compact embedded hypersurfaces $M_t$ satisfying $H>0$ during
the evolution the inequality
\begin{equation*}
2\tau\int_{M_t}\left(\frac{\partial H}{\partial t}-2\,\nabla^{M}H
\cdot \nabla^{M}f
+A(\nabla^{M}f,\nabla^{M}f)-\frac{H}{2\tau}\right) u\,dS \ge 0
\end{equation*}
holds and therefore
\begin{equation*}
\frac{d}{dt}\mathcal{W}_H(\Omega_t, f(t), \tau(t))
\ge 2\tau \int_{\Omega_t}\left|\nabla_i\nabla_j f-\frac{\delta_{ij}}{2\tau}\,\right|^2u\,dx\nonumber\\
\end{equation*}
for $\tau = a-t$ where $a\ge T$ and $t<T$. In particular, this
leads to the above lower volume ratio bound in this case. }

\medskip

Note that the expression
\begin{equation*}
Z(\nabla^Mf)\equiv\frac{\partial H}{\partial
t}-2\,\nabla^{M}H\cdot\nabla^{M}f+A(\nabla^{M}f,\nabla^{M}f)
\end{equation*}
is the central quantity in Hamilton's Harnack inequality for
convex solutions of the mean curvature flow (see \cite{Ha}). Even
though $Z(\nabla^Mf)$ vanishes on translating solutions for
$u=e^{x_{n+1}-\tau}$ our calculations will, due to the
non-compactness of $\Omega_t$ and the non-integrability of all
integrands in this case, not lead to $Z(\nabla^Mf)-H/2\tau$ on the
right hand side.

A direct calculation shows that regions bounded by certain eternal
solutions of mean curvature flow, such as the product of
$\mathbb{R}^{n-1}$ with the grim reaper curve given by
$y=-\log\cos x +t$, do not satisfy the lower volume bound
statement for large $r$ and hence, should the conjecture hold,
cannot occur as a rescaling limit in this situation. Similarly,
certain stationary (zero mean curvature) hypersurfaces would then
be ruled out as rescaling limits such as for instance the catenoid
minimal surface in $\mathbb{R}^3$ and two parallel hyperplanes. In
the positive mean curvature case, White (\cite{Wh}) has previously
shown that certain solutions of mean curvature flow, in particular
the grim reaper hypersurface, cannot occur as rescaling limits.

The embeddedness assumption for the hypersurfaces $M_t$ is essential. Angenent (\cite{A}) has shown, that solutions of the curve-shortening flow
with self-intersections have the grim reaper curve as rescaling limit near singularities.

This paper is organised as follows. In Section 2, we define
entropies for open subsets $\Omega$ of complete (possibly
non-compact) Riemannian manifolds with respect to a given smooth
function $\beta$ defined on $\partial\Omega$ and establish some of
their properties.

In Section 3, we derive the entropy formula involving the Harnack
expression for evolving domains in $\rn$. All of the calculations
go through with necessary modifications such as adding Ricci and
scalar curvature terms in the appropriate places in the case of a
fixed ambient manifold or a background Ricci flow solution.
However, at the moment we do not see how they might lead to
equally interesting consequences.

In Section 4, we state our conjecture and show several
consequences it would lead to, such as a lower local volume ratio
bound and non-existence of certain degenerate rescaling limits.

In Appendix A, we give some explicit examples of entropy
functionals and values in $\rn$.

In the paper, a version of the logarithmic Sobolev inequality on bounded open sets $\Omega$ in complete Riemannian manifolds is used. In Appendix B, we
provide a proof based on the standard Sobolev inequality, essentially following Gross (\cite{G}).

In Appendix C, we give a derivation of a Harnack type evolution equation associated with solutions of a backward heat equation. This equation is one of the
central results in \cite{P} and is also one of the main ingredients in the proof of our entropy formula. Details of this calculation first appeared in \cite{KL}
and \cite{N}.

The work presented in this paper was inspired by a discussion with
Grisha Perelman in January 2003 in Berlin. I would like to thank
Richard Hamilton, Gerhard Huisken, Dan Knopf, Oliver Schn\"urer,
Carlo Sinestrari, Peter Topping, Mu-Tao Wang and Brian White for
helpful discussions. I am particularly indebted to Felix Schulze
for a number of valuable suggestions.

\bigskip

\section{Entropy type functionals for domains in Riemannian manifolds}\label{entropy}

For open subsets $\Omega$ of an $(n+1)$ - dimensional complete (possibly non-compact) Riemannian
manifold $(X, g)$, functions $f:\bar\Omega\to\mathbb{R}$ and
$\beta:\partial\Omega\to\mathbb{R}$ and $\tau >0$ we consider the
quantity
\begin{equation*}
\mathcal{W}_\beta(\Omega, g, f, \tau)=\int_\Omega\left(\tau(|\nabla f|^2+R)+f-(n+1)\right)\,u\,dV+2\tau\int_{\partial\Omega}\beta u \,dS
\end{equation*}
where
\begin{equation*}
u=\frac{e^{-f}}{(4\pi\tau)^{\frac{n+1}{2}}}.
\end{equation*}
The scalar curvature $R$, the expression $|\nabla f|^2$ and the
volume and area elements $dV$ and $dS$ are taken with respect to
the metric $g$. We then define an associated entropy by
\begin{equation*}
\mu_\beta (\Omega, g, \tau)
 = \inf \left\{
\mathcal{W}_\beta (\Omega, g, f, \tau)\, ,\int_\Omega u\,dV =1 \right\}.
\end{equation*}

For $\beta = 0$ and $\Omega = X$, $\mathcal{W}_\beta(\Omega, g, f, \tau)$  and $\mu_\beta (\Omega, g, \tau)$ reduce to Perelman's functional $\mathcal{W}(g, f,
\tau)$ and his entropy quantity $\mu(g, \tau)$. We therefore write $\mathcal{W}$ for $\mathcal{W}_0$ and $\mu$ for $\mu_0$. We use $n+1$ instead of $n$ as we
will later be interested mainly in the hypersurface $\partial\Omega$ which we prefer to be $n$-dimensional.

When we do not intend to vary the metric we consider
\begin{equation*}
\mathcal{W}_\beta (\Omega, f, \tau)=\int_\Omega\left(\tau |\nabla f|^2+f-(n+1)\right)\,u\,dV+2\tau\int_{\partial\Omega}\beta u \,dS
\end{equation*}
with infimum $\mu_\beta (\Omega, \tau)$.

We shall only consider sets with smooth boundaries and smooth
functions $f$ and $\beta$ although the above expressions also make
sense for more general sets and functions. In case $\Omega$ is
unbounded we require suitable integrability conditions on $f$ and
$\beta$. The function $\beta$ could be the restriction to
$\partial\Omega$ of a function on $X$ or be defined only on
$\partial\Omega$. An important example of the latter is $\beta =
H$ where $H$ is the mean curvature of $\partial\Omega$ with
respect to the outer unit normal.

In this section, we derive several basic properties for these entropies. Some specific examples including calculations of entropy values for some natural choices
of sets in $\rn$ are discussed in Appendix A.

\begin{pro}\label{pro:lowerbound}
Suppose that $\Omega$ is bounded with smooth boundary and that $\beta$ is smooth. Then for any $\tau >0$ we have
\begin{equation*}
\mu_\beta (\Omega, g, \tau)\ge -c(n, \Omega, g)\left(1+ \log (1+\tau)+ \tau\,\sup_{\partial\Omega}|\beta|(1+\sup_{\partial\Omega}|\beta|)\right).
\end{equation*}
The same lower bound holds for $\mu_\beta (\Omega, \tau)$.
\end{pro}

\begin{rem}\label{rem:c}
{\rm \noindent  The lower bound for $\mu_\beta(\Omega, g, \tau)$ and for $\mu_\beta(\Omega, \tau)$ follows from the logarithmic Sobolev inquality for $\Omega$
which in turn can be derived from the standard Sobolev inequality (see Appendix B). The constant $c(n, \Omega, g)$ thus depends on the constant in the Sobolev
inequality and the $L^1(\partial\Omega)$ - trace inequality for $C^1(\bar\Omega)$ - functions, the latter controlling the boundary integral. The metric enters
via bounds for the Riemann curvature tensor on $\bar\Omega$ and the explicit bound for the $\sup_\Omega |R|$ - term arising from the functional. The Proposition
holds for more general sets such as bounded sets of finite perimeter and for bounded $\beta$. }
\end{rem}

\noindent{\bf Proof of Proposition \ref{pro:lowerbound}.} We give the proof only for $\mu_\beta(\Omega, g, \tau)$. For $\mu_\beta(\Omega, \tau)$ simply set the
scalar curvature term to zero. We essentially modify the arguments in \cite{KL} and \cite{N}.

Setting $u=\varphi^2$ and using the condition $\int_\Omega u\,dV =1$ we obtain
\begin{equation}\label{eq:transform}\begin{split}
\mathcal{W}_\beta(\Omega, g, f, \tau)
&=\int_\Omega \left(\tau(4|\nabla\varphi|^2+R\varphi^2)-\varphi^2\log\varphi^2\right)\,dV\\
&\qquad\qquad\qquad\qquad+2\tau\int_{\partial\Omega}\beta\varphi^2\,dS - c(n)(1+\log\tau)
\end{split}\end{equation}
with $\int_\Omega \varphi^2\,dV =1$. The trace inequality
\begin{equation*}
\int_{\partial\Omega}\varphi^2\,dS\le c_2 \int_\Omega \left(|\nabla\varphi^2|+\varphi^2\right)\,dV
\end{equation*}
with $c_2=c_2(\Omega, g)$ in combination with Young's inequality yields
\begin{equation*}
\left|2\tau\int_{\partial\Omega}\beta\varphi^2\,dS\right|\le \int_\Omega
2\tau|\nabla\varphi|^2\,dV+c_3\,\tau\,\sup_{\partial\Omega}|\beta|\,(1+\sup_{\partial\Omega}|\beta|)
\end{equation*}
where $c_3$ depends on $c_2$.  Here we have used again the condition $\int_\Omega \varphi^2\,dV =1$. Combining this with (\ref{eq:transform}) yields

\begin{equation}\label{eq:x}\begin{split}
\mathcal{W}_\beta(\Omega, g, f, \tau)\ge \int_\Omega &\left(2\tau|\nabla\varphi|^2-\varphi^2\log\varphi^2\right)\,dV\\ &\qquad -
c_4\left(1+\log\tau+\tau\,(\sup_\Omega |R|+\sup_{\partial\Omega}|\beta|(1+\sup_{\partial\Omega}|\beta|))\right)
\end{split}\end{equation}
where $c_4$ depends on the previous constants. Scaling the metric gives
\begin{equation}\label{eq:changevar}
\int_\Omega \left(2\tau|\nabla\varphi|^2-\varphi^2\log\varphi^2\right)\,dV =
\int_\Omega\left(|\nabla\varphi_\tau|_\tau^2-\varphi_\tau^2\log\varphi_\tau^2\right)\,dV_\tau - c(n)(1+\log\tau)
\end{equation}
and
\begin{equation*}
\int_\Omega\varphi_\tau^2\,dV_\tau =1
\end{equation*}
where $\varphi_\tau = (2\tau)^{\frac{n+1}{4}}\varphi$ and $dV_\tau$ and $|\nabla\varphi_\tau|_\tau^2$ are taken with respect to $g_\tau = (2\tau)^{-1} g$.

By scaling the standard Sobolev inequality
\begin{equation*}
\left(\int_\Omega |\psi|^{\frac{n+1}{n}}\,dV\right)^{\frac{n}{n+1}}\le c_S(\Omega, g)\int_\Omega\left(|\nabla\psi|+|\psi|\right)\,dV
\end{equation*}
we see that the Sobolev constant $c_S(\Omega, g_\tau)$ can be estimated by $c_S(\Omega, g)(1+\sqrt\tau)$. Therefore, by the logarithmic Sobolev inequality
applied in $\Omega$ with respect to the metric $g_\tau$ (see Appendix B)
\begin{equation*}
\int_\Omega\left(|\nabla\varphi_\tau|_\tau^2 -\varphi_\tau^2\log\varphi_\tau^2\right)\,dV_\tau\ge -c(n)\left(1+\log c_S(\Omega, g) +\log (1+\tau )\right).
\end{equation*}
Combining this inequality with ({\ref{eq:x}) and (\ref{eq:changevar}), we arrive at
\begin{equation*}
\mathcal{W}_\beta(\Omega, g, f, \tau) \ge - c_5 \left(1+\log (1+\tau)+\tau\,\sup_{\partial\Omega}|\beta|(1+\sup_{\partial\Omega}|\beta|)\right)
\end{equation*}
with $c_5 = c_5(n, \Omega, g)$ and for $f$ satisfying $\int_\Omega u\,dV =1$. This gives the desired lower bound for $\mu_\beta (\Omega, g, \tau)$. \qed

\begin{pro} Let $\Omega$ be bounded with smooth boundary and assume $\beta$ to be smooth. Then for every $\tau > 0$ there exists a unique smooth minimizer for
$\mu_\beta(\Omega, g, \tau)$ and $\mu_\beta(\Omega, \tau)$. The
minimizer depends smoothly on $\Omega, g, \beta$ and $\tau$.
\end{pro}

\noindent{\bf Proof.} We only consider $\mu_\beta(\Omega, g,
\tau)$ again. The argument is analogous as in \cite{FIN}. The
necessary semicontinuity and coercivity in $W^{1,2}(\Omega)$ for
the transformed functional
\begin{equation*}
\mathcal{E}(\varphi)=\int_\Omega \left(\tau (4|\nabla\varphi|^2 +R\varphi^2)-\varphi^2\log\varphi^2\right)\,dV +2\tau\int_{\partial\Omega}\beta\varphi^2\,dS -
c(n)(1+\log\tau)
\end{equation*}
for $u=\varphi^2$ subject to the condition $\int_\Omega\varphi^2\,dV=1$ follow from similar arguments as in the proof of the lower bound for $\mu_\beta(\Omega,
g, \tau)$ given above. The uniqueness and smooth dependence on the data is standard.\qed

\bigskip
\noindent The quantity
\begin{equation*}
W= W(f) = \tau (2\Delta f -|\nabla f|^2+R) +f -(n+1)
\end{equation*}
featured in Ch.9 of \cite{P} and in \cite{N}. It arises naturally
in the Euler-Lagrange equation for the functional
$\mathcal{W}_\beta(\Omega, g, f, \tau)$.

\begin{pro}\label{pro:euler}
The minimizer $f_{min}$ for the functional $\mathcal{W}_\beta(\Omega, g, f, \tau)$ subject to the constraint $\int_\Omega u\,dV =1$ satisfies the Euler-Lagrange
equation
\begin{equation*}
W(f_{min}) = \mu_\beta(\Omega, g, \tau)
\end{equation*}
in $\Omega$ and the natural boundary condition
\begin{equation*}
\langle\nabla f_{min}, \nu\rangle =\beta
\end{equation*}
on $\partial\Omega$. Here $\langle\cdot,\cdot\rangle$ refers to the metric $g$. For the minimizer of $\mathcal{W}_\beta(\Omega, f, \tau)$ we have instead
\begin{equation*}
W(f_{min}) = \mu_\beta(\Omega, \tau)
\end{equation*}
where
\begin{equation*}
W(f) = \tau (2\Delta f -|\nabla f|^2) +f -(n+1).
\end{equation*}
\end{pro}

\noindent {\bf Proof.} Standard computation using Lagrange multipliers.\qed

\bigskip
\begin{rem}
{\rm \noindent The Euler-Lagrange equation for the transformed functional
\begin{equation*}
\mathcal{E}(\varphi)=\int_\Omega \left(\tau (4|\nabla\varphi|^2 +R\varphi^2)-\varphi^2\log\varphi^2\right)\,dV +2\tau\int_{\partial\Omega}\beta\varphi^2\,dS -
c(n)(1+\log\tau)
\end{equation*}
for $\varphi^2 = u$ subject to the condition $\int_\Omega\varphi^2\,dV=1$ is
\begin{equation*}
-4\tau\,\Delta \varphi -2\varphi\log \varphi +\tau R\varphi= \mu(\Omega, g, \tau)+ (n+1)\left(1+\frac{1}{2}\log {(4\pi\tau)}\right)\varphi
\end{equation*}
in $\Omega$ with boundary condition $2\langle\nabla \varphi, \nu\rangle=-\beta\varphi$ on $\partial\Omega$. }
\end{rem}

\begin{pro}\label{pro:integrationbyparts} For any function $f:\bar\Omega\to\mathbb{R}$ satisfying
\begin{equation*}
\langle\nabla f,\nu\rangle =\beta
\end{equation*}
on $\partial\Omega$ with respect to the outer unit normal $\nu$ we have
\begin{equation*}
\mathcal{W}_\beta(\Omega, g, f, \tau)=\int_\Omega Wu\,dV
\end{equation*}
with $W = W(f) = \tau (2\Delta f -|\nabla f|^2+R) +f -(n+1)$ and
\begin{equation*}
\mathcal{W}_\beta(\Omega, f, \tau)=\int_\Omega Wu\,dV
\end{equation*}
for $W = W(f) = \tau (2\Delta f -|\nabla f|^2) +f -(n+1).$
\end{pro}

\noindent {\bf Proof.} The boundary condition implies $\langle\nabla u, \nu\rangle=-\beta u$ on $\partial\Omega$ and hence
\begin{equation*}
\mathcal{W}_\beta(\Omega, g, f, \tau)=\int_\Omega\left(\tau (|\nabla f|^2+R)+f-(n+1)\right)\,u\,dV-2\tau\int_{\Omega}\Delta u \,dV
\end{equation*}
by the divergence theorem. Since
\begin{equation*}
\Delta u = u ( |\nabla f|^2 - \Delta f).
\end{equation*}
the claim follows.\qed

\bigskip
\noindent For the next statement we do not require $\Omega$ to be bounded.

\begin{pro}\label{pro:upperbound}
Suppose that $\mu_\beta(\Omega, g, r^2)\ge -c_0$ or $\mu_\beta(\Omega, r^2)\ge -c_0$. Let $B_r(x_0)\subset (X, g)$ satisfy $V(\Omega\cap B_{r/2}(x_0))>0$,
\begin{equation*}
\frac {V(\Omega\cap B_r(x_0))+r^2\int_{\partial\Omega\cap B_r(x_0)}|\beta|\,dS}{V(\Omega\cap B_{r/2}(x_0))}\le c_1
\end{equation*}
and $r^2 |Rm|\le c_2$ in $\Omega\cap B_{r}(x_0)$ for the Riemann tensor of $g$. Then
\begin{equation*}
\frac{V(\Omega\cap B_r(x_0))}{r^{n+1}}\ge\kappa >0
\end{equation*}
with $\kappa =\kappa (n, c_0, c_1, c_2)$.
\end{pro}

\begin{rem}\label{rem:ub}
\noindent{\rm We will actually prove that  $\mu_\beta(\Omega, g,  r^2)$ and $\mu_\beta(\Omega, r^2)$ are bounded from above by the expression
\begin{equation*}
\log\,\frac{V(\Omega\cap B_r(x_0))}{r^{n+1}}+c\,\frac {V(\Omega\cap B_r(x_0))+r^2\int_{\partial\Omega\cap B_r(x_0)}|\beta|\,dS}{V(\Omega\cap B_{r/2}(x_0))}
\end{equation*}
with $c = c(n, r^2|Rm|)$. From this the claim follows immediately.}
\end{rem}

\noindent{\bf Proof of Proposition \ref{pro:upperbound}.} In the case $\Omega =X$ and $\beta =0$ the proof is sketched in Ch.3 of \cite{P} (see \cite{KL} and
\cite{N} for more details). We proceed along similar lines.

If we set $e^{-f}=a\zeta$ the normalisation condition for $f$ becomes
\begin{equation*}
a=\frac{(4\pi r^2)^{\frac{n+1}{2}}}{\int_\Omega\zeta\,dV}.
\end{equation*}
The functional $\mathcal{W}_\beta(\Omega, f, r^2)$ can then be expressed as
\begin {equation*}
\frac{a}{(4\pi r^2)^{\frac{n+1}{2}}}\int_{\Omega}\left( 4 r^2 \frac{|\nabla\zeta|^2}{\zeta}+r^2R\,\zeta-\zeta\log (a\zeta)\right)\,dV - (n+1) +
2r^2\frac{\int_{\partial\Omega}\beta\zeta\,dS}{\int_\Omega\zeta\,dV}.
\end{equation*}
By approximation, we may substitute functions $\zeta\in C^2_0(X)$ into this expression. We choose as $\zeta$ a cut-off function for $B_{r/2}(x_0)$ that is
$\zeta$ satisfies $\chi_{B_{r/2}(x_0)}\le\zeta\le\chi_{B_r(x_0)}$ as well as
\begin{equation*}
4r^2\frac{|\nabla\zeta|^2}{\zeta}\le 8 r^2 \sup |\nabla^2\zeta|\le
c
\end{equation*}
where $c$ is a constant which depends on $r^2\sup_{ \Omega\cap B_r(x_0)} |Rm|$ and is therefore bounded by $c_2$. Since
\begin{equation*}
\int_\Omega\zeta\,dV \ge V(\Omega\cap B_{r/2}(x_0))>0
\end{equation*}
we can thus estimate
\begin {equation*}
\frac{1}{(4\pi r^2)^{\frac{n+1}{2}}}\int_{\Omega}4 r^2 a \frac{|\nabla\zeta|^2}{\zeta}\,dV\le c\frac{V(\Omega\cap {\rm spt}\;\zeta )}{\int_\Omega\zeta\,dV}\le
c\frac{V(\Omega\cap B_r(x_0))}{V(\Omega\cap B_{r/2}(x_0))}.
\end{equation*}
Jensen's inequality now implies
\begin{equation*}
-\frac{1}{(4\pi r^2)^{\frac{n+1}{2}}}\int_\Omega a\zeta\log (a\zeta)\,dV \le -\frac{1}{(4\pi r^2)^{\frac{n+1}{2}}}\int_\Omega a\zeta\,dV \log
\left(\frac{1}{V(\Omega\cap {\rm spt} \zeta)}\int_\Omega a\zeta\,dV\right).
\end{equation*}
Since ${\rm spt}\;\zeta =\overline{B_r(x_0)}$ and in view of the normalisation condition the right hand side equals
\begin{equation*}
\log\left(\frac{V(\Omega\cap B_r(x_0))}{(4\pi r^2)^{\frac{n+1}{2}}}\right).
\end{equation*}
The scalar curvature integral is estimated using the boundedness assumption on the Riemann tensor in $\Omega\cap B_r(x_0)$. This yields the upper bound for $\mu
(\Omega, g, r^2)$ and $\mu (\Omega, r^2)$ stated in Remark \ref{rem:ub}. \qed

\begin{rem}\noindent{\rm In \cite{P}, Perelman ruled out the occurrence of collapsed metrics as rescaling limits of compact, finite time solutions of Ricci flow.
A metric $g$ on $X$ is called \emph{collapsed} if there exists a sequence of balls
$B_{r_k}(x_k)\subset (X, g)$ satisfying $r_k^2 |Rm|\le 1$ in $B_{r_k}(x_k)$ for which
\begin{equation*}
\frac{V(B_{r_k}(x_k))}{r_k^{n+1}}\to 0.
\end{equation*}
An important example of a collapsed metric is the so-called \emph{cigar soliton} solution of the Ricci flow given by $X=\mathbb{R}^2$ endowed with the metric
\begin{equation*}
ds^2=\frac{dx^2+dy^2}{1+x^2+y^2}.
\end{equation*}
On collapsed metrics we have $\inf_{\tau >0}\mu_\beta(g, \tau) = -\infty$ by the proposition. }\end{rem}

The following reformulation of Proposition \ref{pro:upperbound}
links a kind of volume collapsing behaviour of subsets of $(X, g)$
to a property of the entropy $\mu_\beta(\Omega, g, \tau)$.

\begin{cor}\label{cor:-infty}
{\noindent If for some fixed constants $c_1$ and $c_2$ we can find a sequence of balls $B_{r_k}(x_k)\subset (X, g)$ such that $V(\Omega\cap B_{r_k/2}(x_k))>0$,
\begin{equation*}
\frac {V(\Omega\cap B_{r_k}(x_k))+r_k^2\int_{\partial\Omega\cap B_{r_k}(x_k)}|\beta|\,dS}{V(\Omega\cap B_{r_k/2}(x_k))}\le c_1,
\end{equation*}
$r_k^2 |Rm|\le c_2$ in $\Omega\cap B_{r_k}(x_k)$ and
\begin{equation*}
\frac{V(\Omega\cap B_{r_k}(x_k))}{r_k^{n+1}}\to 0
\end{equation*} then
$\inf_{\tau >0}\mu_\beta(\Omega, g, \tau) = -\infty$ and $\inf_{\tau >0}\mu_\beta(\Omega, \tau) = -\infty$. }
\end{cor}

For compact $\Omega$ we can of course always find such a sequence
of balls with radii tending to infinity. In the case of
non-compact regions the sitation is more interesting. Examples are the
following regions in $X=\rn$:

\medskip

\noindent (1) The slab
\begin{equation*}
\Omega = \{x\in\rn,\,-d<x_{n+1}<d\}
\end{equation*}
for some $d>0$. On the hypersurface $M=\partial\Omega$ we have
$H=0$. The enclosed region $\Omega$ satisfies $V(\Omega\cap
B_{r/2})>0$ and
\begin{equation*}
\frac{V(\Omega\cap B_r)}{V(\Omega\cap B_{r/2})}\le c(n, d)
\end{equation*}
for all balls $B_r=B_r(0)$. Moreover,
\begin{equation*}
\lim_{r\to\infty}\frac{V(\Omega\cap B_r)}{r^{n+1}} = 0.
\end{equation*}

\noindent (2) The 'smaller' of the two regions bounded by the
catenoid minimal surface $M=\partial\Omega$ in $\mathbb{R}^3$
given by
\begin{equation*}
\Omega =\{x=(\hat x, x_3)\in\mathbb{R}^3, \,|\hat x|\ge
1,\,|x_3|\le \cosh^{-1} |\hat x|\}.
\end{equation*}

Note that $H=0$ on $\partial\Omega$. One checks that there is a
constant $c_1$ such that for all $r\ge 2$
\begin{equation*}
\frac{V(\Omega\cap B_r)}{V(\Omega\cap B_{r/2})}\le c_1
\end{equation*}
and
\begin{equation*}
V(\Omega\cap B_r)\le c_1 r^2\log (1+r)
\end{equation*}
so that
\begin{equation*}
\lim_{r\to\infty}\frac{V(\Omega\cap B_r)}{r^3} = 0.
\end{equation*}

\noindent (3) The translating solution of mean curvature flow
corresponding to the grim reaper hypersurface $M=\partial\Omega$
where $\Omega=\mathbb{R}^{n-1} \times G$ with
\begin{equation*}
G=\left\{(x_n, x_{n+1})\in\mathbb{R}^2,\, -\pi/2 < x_n < \pi/2,
\,x_{n+1}>-\log\cos  x_n\right\}.
\end{equation*}

An explicit calculation shows that the mean curvature satisfies
$H(x) = e^{-x_{n+1}}$ for any $x\in M=\partial\Omega$ . One
therefore checks directly that there is a sequence of balls
$B_{r_k}(x_k)$ with $r_k\to\infty$ satisfying $V(\Omega\cap
B_{r_k/2}(x_k))>0$,
\begin{equation*}
\frac{V(\Omega\cap B_{r_k}(x_k))}{V(\Omega\cap B_{r_k/2}(x_k))}\le
c(n),
\end{equation*}
\begin{equation*}
\frac{r_k^2\int_{\partial\Omega\cap
B_{r_k}(x_k)}H\,dS}{V(\Omega\cap B_{r_k/2}(x_k))}\le 1
\end{equation*}
and
\begin{equation*}
\frac{V(\Omega\cap B_{r_k}(x_k))}{r_k^{n+1}}\to 0.
\end{equation*}

\bigskip

\section{An Entropy type formula for evolving domains in $\rn$}\label{ricci}

In this section we restrict ourselves to domain evolution in
$\rn$. All the calculations go through for fixed Riemannian
manifolds or Ricci flow solutions as ambient space if we add Ricci
and scalar curvature terms in the appropriate places. Howerer, in
this case the formulas do not immediately seem to lead to any
interesting consequences.

We evolve bounded open subsets $(\Omega_t)_{t\in [0, T)}$ with
smooth boundary hypersurfaces $(M_t)_{t\in [0, T)}$ in $\rn$. More
precisely, $\bar\Omega_t=\phi_t(\bar\Omega)$ with
$M_t=\partial\Omega_t=\phi_t(\partial\Omega)$ where $\phi_t=\phi
(\cdot, t):\bar \Omega\to \rn \,, t\in [0, T)$ is a smooth
one-parameter family of diffeomorphisms. We will often abbreviate
\begin{equation*}
x=\phi (p, t)
\end{equation*}
for $p\in\bar\Omega$. The normal speed of $M_t$ with respect to the inward pointing normal $-\nu$ is defined by
\begin{equation*}
\beta=\beta_{M_t}=-\frac{\partial x}{\partial t}\cdot \nu
\end{equation*}
for $x\in M_t$ or expressed in terms of the
embedding map $\phi (\cdot, t)$ by
\begin{equation*}
\beta(p, t)=-\frac{\partial \phi}{\partial t}(p, t)\cdot
\nu(\phi(p, t))
\end{equation*}
for $p\in\partial\Omega$. We assume the function $\beta$ to be
smooth. If for instance $\beta=H$, the mean curvature of $M_t$,
this describes mean curvature flow up to diffeomorphisms
tangential to $M_t$.

Let us assume more specifically that the family of subsets
$(\Omega_t)_{t\in (0, T)}$ evolves by the equation
\begin{equation}\label{eq:omegaflow}
\frac{\partial x}{\partial t}=-\nabla f(x, t)
\end{equation}
for $x\in\Omega_t$. This flow is compatible with the evolution of
the boundaries $M_t=\partial \Omega_t$ with normal speed $\beta$
if $f$ satisfies the condition $\nabla f\cdot\nu =\beta$ on $M_t$.
Suppose $f(t)$ satisfies the equation
\begin{equation}\label{eq:fequation}
\left(\frac{\partial }{\partial t} +\Delta\right) f = |\nabla
f|^2+\frac{n+1}{2\tau}
\end{equation}
in $\Omega_t$ for $t\in (0, T)$. The total time derivative  of $f$ is given by
\begin{equation}\label{eq:totalder}
\frac{df}{dt}=\frac{\partial f}{\partial t}+\nabla f\cdot
\frac{\partial x}{\partial t}=\frac{\partial f}{\partial
t}-|\nabla f|^2.
\end{equation}
Hence (\ref{eq:fequation}) can also be written as
\begin{equation}\label{eq:totalfequation}
\left(\frac{d}{dt} +\Delta\right) f = \frac{n+1}{2\tau}.
\end{equation}
If $\tau(t)>0$ evolves by $\frac{\partial\tau}{\partial t}=-1$ then (\ref{eq:fequation}) is equivalent to the equation
\begin{equation}\label{eq:uequation}
\left(\frac{\partial }{\partial t}+\Delta\right)u=0
\end{equation}
for
\begin{equation*}
u=\frac{e^{-f}}{(4\pi\tau)^{\frac{n+1}{2}}}.
\end{equation*}

The above equations are more precisely expressed in terms of the
pull back of the function $f$ via the diffeomorphisms evolving
$\Omega_t$. In fact, if we set $x=\phi (q, t)$ where $\phi_t =\phi
(\cdot, t):\Omega\to\Omega_t$, the pulled back function given by
\begin{equation*}
\tilde f(q, t)=f (\phi(q, t), t)
\end{equation*}
satisfies
\begin{equation*}
\frac{df}{dt}(x, t)=\frac{\partial\tilde f}{\partial t}(q, t).
\end{equation*}

\medskip
Analogously to Ch.9 in \cite{P} (see also \cite{N}) the function
$W=\tau(2\Delta f-|\nabla f|^2)+f-(n+1)$ satisfies a nice
evolution equation:

\begin{pro}\label{pro:wequation} Let $(\Omega_t)_{t\in (0, T)}$ be a family of subsets evolving by
(\ref{eq:omegaflow}) that is according to the negative gradient of
functions $f(t)$ satisfying equation (\ref{eq:fequation}). Suppose
also that $\tau(t)>0$ evolves by $\frac{\partial\tau}{\partial
t}=-1$ for $t\in (0, T)$. Then the function
\begin{equation*}
W=\tau(2\Delta f-|\nabla f|^2)+f-(n+1)
\end{equation*}
satisfies the evolution equation
\begin{equation*}
\left(\frac{d}{dt} +\Delta\right)W= 2\tau\left|\nabla_i\nabla_j
f-\frac{\delta_{ij}}{2\tau}\right|^2 + \nabla W\cdot\nabla f
\end{equation*}
in $\Omega_t$.
\end{pro}

\noindent{\bf Proof.} We use Perelman's identity
\begin{equation*}
\left(\frac{\partial}{\partial t} +\Delta\right)W= 2\tau
\left|\nabla_i\nabla_j f-\frac{\delta_{ij}}{2\tau}\;\right|^2+
2\,\nabla W\cdot\nabla f
\end{equation*}
from Ch.9 in \cite{P}. A derivation of this can be found in
\cite{KL} and in \cite{N}. In our evolving coordinates $x=\phi(q,
t)$ we change to total time derivatives for $W$ via
\begin{equation*}
\frac{dW}{dt}=\frac{\partial W}{\partial t}-\nabla W\cdot\nabla f
\end{equation*}
which yields the result. For the convenience of the reader, we repeat the details of the calculation in \cite{N} for the expression $\left(\frac{d}{dt}
+\Delta\right)W$ on evolving sets $\Omega_t\subset\rn$ in Appendix C .\qed

\begin{pro}\label{pro:dtwbeta}
Suppose the conditions of the previous proposition hold. Then
\begin{equation*}
\frac{d}{dt}\int_{\Omega_t}u\,dx = 0
\end{equation*}
for all $t\in (0, T)$. If $f$ satisfies additionally $\nabla
f\cdot \nu =\beta$ on $M_t=\partial\Omega_t$ then
\begin{equation}\label{eq:a}
\frac{d}{dt}\mathcal{W}_\beta(\Omega_t, f(t), \tau(t)) =2\tau
\int_{\Omega_t}\left|\nabla_i\nabla_j
f-\frac{\delta_{ij}}{2\tau}\,\right|^2u\,dx -\int_{M_t}\nabla
W\cdot\nu \,u\,dS
\end{equation}
where $W=\tau (2\Delta f -|\nabla f|^2) +f -(n+1)$.
\end{pro}

\noindent {\bf Proof.} In view of the family of diffeomorphisms
generated by
\begin{equation*}
\frac{\partial x}{\partial t}=-\nabla f =\frac{1}{u}\nabla u
\end{equation*}
the volume element $dx$ on the evolving sets $\Omega_t$ changes by
\begin{equation*}
\frac{d}{dt}\,dx = -\Delta f\,dx.
\end{equation*}
Since also
\begin{equation*}
\frac{du}{dt}=\frac{\partial u}{\partial t}+\frac{|\nabla u|^2}{u}
\end{equation*}
and
\begin{equation*}
\Delta u = (|\nabla f|^2-\Delta f)\,u
\end{equation*}
we obtain in $\Omega_t$
\begin{equation}\label{eq:dtudx}
\frac{d}{dt}\left(u\,dx\right)=\left(\frac{\partial u}{\partial
t}+\Delta u\right)dx = 0
\end{equation}
by equation (\ref{eq:uequation}). Thus
\begin{equation*}
\frac{d}{dt}\int_{\Omega_t}u\,dx = 0.
\end{equation*}

Combining the Neumann boundary condition, Proposition \ref{pro:integrationbyparts}, identity (\ref{eq:dtudx}) and the evolution equation for $W$ in Proposition
\ref{pro:wequation} we then calculate
\begin{align}
&\frac{d}{dt}\mathcal{W}_\beta(\Omega_t, f(t), \tau(t))\nonumber\\
&\qquad\qquad= \frac{d}{dt}\int_{\Omega_t}Wu\,dx
= \int_{\Omega_t}\left(\frac{d}{dt} +\Delta\right)W\, u\,dx -\int_{\Omega_t}\Delta W \,u\,dx\nonumber\\
&\qquad\qquad= \int_{\Omega_t}2\tau \left|\nabla_i\nabla_j
f-\frac{\delta_{ij}}{2\tau}\right|^2\,u\,dx-
\int_{\Omega_t}\left(\nabla W\cdot\nabla u +\Delta W \,u
\right)\,dx\nonumber
\end{align}
where we again used $\nabla u = -u \nabla f$. The last integral
equals
\begin{equation*}
-\int_{\Omega_t}{\textrm div}\left(\nabla W u\right)\,dx.
\end{equation*}
The result then follows by applying the divergence theorem. \qed

\begin{rem}
\noindent{\rm For a fixed domain $\Omega$ (that is when $\beta
=0$) inside a Riemannian manifold of nonnegative Ricci curvature
the inequality
\begin{equation*}
\frac{d}{dt}\mathcal{W}(\Omega, f(t), \tau(t)) \ge
-\int_{M_t}\nabla W\cdot\nu \,u\,dS
\end{equation*}
for a solution $f$ of the above backward heat equation appears in
\cite{N}. Ni then shows that
\begin{equation*}
-\langle\nabla W, \nu\rangle = 2\tau A(\nabla^Mf, \nabla^M f)
\end{equation*}
and is therefore non-negative for a convex boundary (see below for
a generalisation of the corresponding calculation to evolving
domains), thus obtaining
\begin{equation*}
\frac{d}{dt}\mathcal{W}(\Omega, f(t), \tau(t))\ge 0.
\end{equation*}}
\end{rem}

\medskip

When examining the integrand $-\nabla W\cdot\nu$ of the above
boundary integral more closely, an interesting relation with the
expression in Hamilton's Harnack inequality for the mean curvature
of a hypersurface evolving by mean curvature flow emerges. To
appreciate this one should first note that the hypersurfaces $M_t$
evolve by the equation
\begin{equation}\label{eq:mtequation}
\frac{\partial x}{\partial t}=-\beta\nu-\nabla^Mf
\end{equation}
due to the Neumann boundary condition for $f$ where $\nabla^M$
denotes the tangential gradient on the hypersurfaces $M_t$.

\begin{pro}\label{pro:harnack}
Under the above conditions on $(M_t)$ and $f(t)$ the quantity $W$
satisfies the identity
\begin{equation}\label{eq:gradw}
-\nabla W\cdot\nu=2\tau\left(\frac{\partial\beta}{\partial
t}-2\,\nabla^{M}\beta\cdot\nabla^{M}f+A(\nabla^{M}f,\nabla^{M}f)-\frac{\beta}{2\tau}\right)
\end{equation}
for all $t<T$, $a\ge T$ and $\tau=a-t$ where $A$ denotes the
second fundamental form of $M_t$. This implies the inequality
\begin{equation*}
\frac{d}{dt}\mathcal{W}_\beta (\Omega_t, f(t), \tau(t)) \ge
2\tau\int_{M_t}\left(\frac{\partial\beta}{\partial
t}-2\,\nabla^{M}\beta \cdot \nabla^{M}f
+A(\nabla^{M}f,\nabla^{M}f)-\frac{\beta}{2\tau}\right)u\,dS
\end{equation*}
\end{pro}

\noindent{\bf Proof.} In view of equation
(\ref{eq:totalfequation}) we have
\begin{equation*}
W=-\tau\left(2\frac{df}{dt}+|\nabla f|^2\right)+f.
\end{equation*}
We now calculate similarly as in Appendix C
\begin{equation*}
\frac{d}{dt}\nabla f =\nabla^2f\left(\nabla f,
\,\cdot\,\right)+\nabla \frac{df}{dt}.
\end{equation*}
A calculation as for instance in (\cite{Hu1}) using the evolution
equation (\ref{eq:mtequation}) for the hypersurfaces $M_t$ yields
\begin{equation*}
\frac{d\nu}{dt}=\nabla^{M}\beta-A(\nabla^{M}f,\,\cdot\,)
\end{equation*}
for the outward unit normal field on $M_t$. The second term arises
from the definition of $A$ in terms of tangential derivatives of
$\nu$. Combining these and differentiating the identity $\beta =
\nabla f\cdot\nu$ yields
\begin{equation*}
\frac{d\beta}{dt} =\nabla^2 f\left(\nabla f, \nu\right)+\nabla
\frac{df}{dt}\cdot\nu+\nabla^{M}\beta\cdot\nabla^{M}f-A(\nabla^{M}f,\nabla^{M}f).
\end{equation*}
Since
\begin{equation*}
\nabla W\cdot\nu=-\tau\left(2\nabla\frac{df}{dt}\cdot\nu+\nabla
|\nabla f|^2\cdot\nu\right)+\nabla f\cdot\nu
\end{equation*}
and $\nabla |\nabla f|^2\cdot\nu= 2\,\nabla^2f\left(\nabla f,
\nu\right)$ we obtain the result by observing
\begin{equation*}
\frac{d\beta}{dt}=\frac{\partial\beta}{\partial
t}-\nabla^{M}\beta\cdot\nabla^{M}f
\end{equation*}
in view of (\ref{eq:mtequation}). The integral inequality then
follows from Proposition \ref{pro:dtwbeta}.\qed

\begin{rem}\noindent{\rm  Let $f^{t_0}$ be the minimizer for $\mu_\beta(\Omega_{t_0},
\tau(t_0))$. Since $W(f^{t_0})\equiv {\rm constant}$ (see
Proposition \ref {pro:euler}) we have
\begin{equation*}
\int_{M_{t_0}}\nabla W\cdot\nu\,u\,dS =0
\end{equation*}
at time $t_0$. However, even if we assume that $f(t)$ for $t<t_0$
satisfies the 'end' condition $f(t_0)=f^{t_0}$ we cannot conclude
that
\begin{equation*}
\lim_{t\to t_0}\int_{M_{t}}\nabla W\cdot\nu\,u\,dS =0
\end{equation*}
and that therefore (note that $\mathcal{W}_\beta$ is
differentiable at $t_0$)
\begin{equation*}
\frac{d}{dt}_{|{t_0}}\mathcal{W}_\beta(\Omega_t, f(t), \tau(t))\ge
0.
\end{equation*}
The problem occurs since $\nabla W$ involves third derivatives of
$f$ which won't behave continuously on the boundary for $t\to t_0$
unless we impose some kind of higher order compatibility condition
on the 'end' data $f^{t_0}$ on $M_{t_0}=\partial\Omega_{t_0}$.}
\end{rem}

\bigskip

\section{A conjectured Harnack type inequality for mean curvature flow and its consequences}

For $\beta = H$, the expression
\begin{equation*}
Z(\nabla^Mf)\equiv\frac{\partial H}{\partial
t}-2\,\nabla^{M}H\cdot\nabla^{M}f+A(\nabla^{M}f,\nabla^{M}f)
\end{equation*}
in Proposition \ref{pro:harnack} is the central quantity in
Hamilton's Harnack inequality for convex solutions of the mean
curvature flow (see \cite{Ha}). Hamilton showed, that $Z(V)$
vanishes on translating solutions of mean curvature flow for some
vector field $V$ which is tangential to the hypersurfaces $M_t$.
His Harnack inequality states that
\begin{equation*}
Z(V)+\frac{H}{2t}\ge 0
\end{equation*}
holds for any tangential vector field $V$ on a convex solution of
mean curvature flow for $t>0$ with equality for a suitable vector
field on a homothetically expanding solution. We observe that on
homothetically shrinking solutions that is where
\begin{equation*}
H=\frac{x\cdot\nu}{2\tau}
\end{equation*}
the identity
\begin{equation*}
2\tau Z(V)- H =0
\end{equation*}
holds for $V=\nabla^M f$ where $f=|x|^2/4\tau$.

Because of the term $-H$ we cannot expect this expression to be
nonnegative for a general solution and for a general $V$. It
certainly is negative on translating solutions for a suitable $V$.
However, it seems reasonable to expect that this quantity or its
integral over $M_t$ has a sign, at least on compact solutions. For
non-compact solutions, our calculations do not lead to the
integral inequality in Proposition \ref{pro:dtwbeta} since the
integral expressions are usually not well-defined in this case as
will see a little later in the case of translating solutions. The
above considerations lead us to the following

\bigskip
{\bf Conjecture.} {\sl Let $(M_t)_{t<T}$ be a family of compact
embedded hypersurfaces evolving by their mean curvature. Assume
$H>0$ during the flow. Let $\tau = a-t$ for fixed $a\ge T$ and all
$t<T$. Then the Harnack type inequality
\begin{equation*}
2\tau\left(\frac{\partial H}{\partial t}-2\,\nabla^{M}H\cdot V
+A(V,V)\right)-H\ge 0
\end{equation*}
holds on $M_t$ for all $t<T$ and all tangential vector fields $V$
with equality on homothetically shrinking solutions for
$V=\nabla^M f$ and $f=|x|^2/4\tau$. If we rescale mean curvature
flow by considering $\tilde x(s)=1/\sqrt {2\tau (t)}\, x(t)$ for
$s=-\log \sqrt {2\tau (t)}$ then the conjectured inequality
becomes
\begin{equation*}
\frac{\partial\tilde H}{\partial s}-2\tilde{\nabla}^{\tilde
M}\tilde {H}\cdot V + \tilde{A}(V, V)\ge 0
\end{equation*}
for all $s>0$. A weaker form of the conjecture which suffices for
the applications we have in mind is that
\begin {equation*}
{\textrm (C)} \qquad 2\tau\int_{M_t}\left(\frac{\partial
H}{\partial t}-2\,\nabla^{M} H
\cdot\nabla^{M}f+A(\nabla^{M}f,\nabla^{M}f)-\frac{H}{2\tau}\right)\,u\,dS
\ge 0
\end{equation*}
where $f$ satisfies
\begin{equation*}
\left(\frac{\partial }{\partial t} +\Delta\right) f= |\nabla
f|^2+\frac{n+1}{2\tau}
\end{equation*}
in $\Omega_t$ for $t<T$ with the
boundary condition $\nabla f\cdot\nu =H$ and the domains evolve by
a family of diffeomorphisms generated by $-\nabla f$. Note that $f$ blows up like $|x|^2/4(T-t)$ and therefore $|\nabla^Mf|$ like $|x|/(T-t)$ for $t\to T$ near a singularity corresponding to a homothetically shrinking
solution.}

\medskip

Let us give two explicit examples of mean curvature flow solutions
which illustrate the situation: First note that the evolution
equation for the hypersurfaces $M_t$ in the case $\beta=H$ is
\begin{equation*}
\frac{\partial x}{\partial t}=-H\nu-\nabla^Mf,
\end{equation*}
which is mean curvature flow up to tangential diffeomorphisms.

If $\Omega_t$ is the interior of a homothetically shrinking
solution of mean curvature flow, that is up to translation in time
\begin{equation*}
\Omega_t = \sqrt{2\tau}\,\Omega_0
\end{equation*}
for $\tau = T-t$, then $f=|x|^2/4(T-t)$ is a solution of equation
(\ref{eq:fequation}). The Neumann boundary condition above becomes
simply
\begin{equation*}
H=\frac{x\cdot\nu}{2\tau}.
\end{equation*}
In this situation,
\begin{equation*}
\frac{\partial H}{\partial
t}-2\,\nabla^{M}H\cdot\nabla^{M}f+A(\nabla^{M}f,\nabla^{M}f)-\frac{H}{2\tau}=0
\end{equation*}
so $\nabla W\cdot\nu =0$.

For translating solutions of mean curvature flow the
quantity $-\nabla W\cdot\nu$ is negative for positive $H$.
However, our rate of change formula for $\mathcal{W}_H$ does not
hold in this case as the entropy calculations are not justified in
this situation:

Indeed, if $\Omega_t$ is the interior of a translating solution of
mean curvature flow, that is up to rotation in $\rn$
\begin{equation*}
\Omega_t=\Omega + t e_{n+1}
\end{equation*}
for some fixed set $\Omega$ and for all $t\in\mathbb{R}$ then
\begin{equation*}
f=-x_{n+1} +\tau -\log (4\pi\tau)^{\frac{n+1}{2}}
\end{equation*}
solves the boundary value problem. The Neumann boundary condition
on $M_t$ in this case becomes $H=-\nu_{n+1}$.

We note that $M_t$ and $\Omega_t$ are necessarily unbounded since
compact solutions cannot exist for all $t\in\mathbb{R}$ by
comparison with spheres shrinking to points in finite time.
Moreover, the function $u$ featuring in the integrand of the
entropy functional as well as in the normalisation condition
required for the entropy is given by
\begin{equation*}
u=\frac{e^{-f}}{(4\pi\tau)^{\frac{n+1}{2}}}=e^{x_{n+1}-\tau}
\end{equation*}
in our example. In view of the comparison principle for mean
curvature flow applied to $M_t$ and hyperplanes
$\{x\in\rn,\,x_{n+1}=a\}$,  which are stationary solutions of mean
curvature flow, the sets $\Omega_t$  have an unbounded
intersection with the upper half space $\{x\in\rn,\,x_{n+1}>0\}$
for every $t\in\mathbb{R}$. Therefore, the function $u$ is an
illegal choice in the normalisation condition
$\int_{\Omega_t}u\,dx =1$  as it is not integrable on $\Omega_t$.

\medskip

There are a number of important consequences of inequality (C)
especially for the open problem of no local volume collapse for
mean curvature flow solutions (an analogue of Perelman's no local
collapsing for Ricci flow solutions) and consequently
non-existence of certain degenerate rescaling limits. This should
provide sufficient motivation for settling the conjecture.

\begin{pro} Suppose that the conjectured inequality (C) holds. Then
\begin{equation*}
\frac{d}{dt}\mathcal{W}_H(\Omega_t, f(t), \tau(t))\ge 0
\end{equation*}
for $t<t_0$ and therefore the entropy is monotonic that is
\begin{equation*}
\mu_H(\Omega_{t_1}, a-t_1)\le \mu_H(\Omega_{t_2}, a-t_2)
\end{equation*}
for $0\le t_1\le t_2 <T$ and any $a\ge T$.
\end{pro}

{\bf Proof.} The first inequality follows directly from
Proposition \ref{pro:harnack} applied to $\beta =H$ and from (C).
To derive the second inequality we let $f^{t_0}$ for $t_0<T$ be
the minimizer for $\mu_H(\Omega_{t_0}, \tau(t_0))$ and let $f(t)$
in addition to the equation
\begin{equation*}
\left(\frac{\partial }{\partial t} +\Delta\right) f= |\nabla
f|^2+\frac{n+1}{2\tau}
\end{equation*} in $\Omega_t$ and the boundary condition $\nabla f\cdot\nu=H$ on $\partial\Omega_t$ for $t<t_0$ satisfy the 'end'
condition $f(t_0)=f^{t_0}$. Since
\begin{equation*}
\frac{d}{dt}\mathcal{W}_H(\Omega_t, f(t), \tau(t))\ge 0
\end{equation*} we have
\begin{equation*}
\mathcal{W}_H(\Omega_t, f(t), \tau(t))\le
\mathcal{W}_H(\Omega_{t_0}, f(t_0),\tau(t_0))=
\mathcal{W}_H(\Omega_{t_0}, f^{t_0},
\tau(t_0))=\mu_H(\Omega_{t_0},\tau(t_0)).
\end{equation*}
Taking the infimum on the left hand side over all functions
satisfying the normalisation condition
\begin{equation*}
\int_{\Omega_t}\frac{e^{-f}}{(4\pi\tau)^{\frac{n+1}{2}}}\,dV=1
\end{equation*}
we obtain the desired inequality for
the entropies at $t$ and at $t_0$. Since $t$ and $t_0$ were arbitrary
we are done. \qed

\bigskip

\begin{cor}\label{cor:mcf} Let $(M_t)_{t\in [0, T)}$ be a solution of mean curvature flow consisting of smooth, compact, embedded hypersurfaces which enclose
bounded regions $(\Omega_t)_{t\in [0, T)}$ in $\rn$. Let $\tau
=a-t$ for arbitrary but fixed $a\ge T$ and all $t<T$. Suppose
furthermore that the Harnack type inequality (C) holds. Then for
every $r>0$ and every $t\in [0, T)$
\begin{equation*}
\mu_{H}(\Omega_t, r^2)\ge \mu_{H}(\Omega_0, t+r^2).
\end{equation*}
Since $T<\infty$ we have for every $t\in [0, T)$ and $r\in (0, \sqrt{T} ]$
\begin{equation*}
\mu_{H}(\Omega_t, r^2)\ge -c_0
\end{equation*}
where $c_0$ depends only on $n, \Omega_0, T$ and $\sup_{M_0}|H|$. In particular, there is a constant $\kappa >0$ depending only on $n, \Omega_0, T,
\sup_{M_0}|H|$ and $c_1$ such that the inequality
\begin{equation*}
\frac{V(\Omega_t\cap B_r(x_0))}{r^{n+1}}\ge\kappa
\end{equation*}
holds for all $t\in [0, T)$ and $r\in (0, \sqrt{T}]$ in balls $B_r(x_0)$ satisfying the conditions $V(\Omega_t\cap B_{r/2}(x_0))>0$ and
\begin{equation*}
\frac{V(\Omega_t\cap B_r(x_0))+r^2\int_{M_t\cap B_r(x_0)}|H|\,dS}{V(\Omega_t\cap B_{r/2}(x_0))}\le c_1.
\end{equation*}
\end{cor}

{\bf Proof.} By the monotonicity of the entropy and Proposition
\ref {pro:lowerbound} applied with $a=r^2+t,\, t_1=0$ and $t_2=t$
we have

\begin{equation*}
\mu_H(\Omega_t, r^2)\ge \mu_H(\Omega_0, t+r^2)\ge -c(n, T, \sup_{M_0}|H|, \Omega_0)
\end{equation*}
for all $r\le \sqrt{T}$ and $t <T$. The lower volume ratio bounds
then follow from Proposition \ref{pro:upperbound} applied to
$\Omega_t$. \qed
\bigskip

\bigskip

\noindent For $\lambda_j\searrow 0, \,t_j\nearrow T$ and
$x_j\in\rn$ we define a sequence $(\Omega^j_s)$ of rescaled flows
\begin{equation*}
\Omega_s^j =\frac{1}{\lambda_j}\left(\Omega_{\lambda_j^2 s+t_j}-x_j\right)
\end{equation*}
where $s\in (-\lambda_j^{-2}t_j, \lambda^{-2}(T-t_j))\equiv (a_j, b_j)$.

\begin{defi}{\rm \noindent Let $(M_t)_{t\in [0, T)}$ be a compact, smooth, embedded solution of mean curvature flow enclosing bounded regions
$(\Omega_t)_{t\in [0, T)}$ in $\rn$. We call a smooth, embedded solution $(M'_s)_{s\in (-\infty, b)}$ of mean curvature flow enclosing (not necessarily bounded)
$(\Omega'_s)_{s\in (-\infty, b)}$
a \emph{rescaling limit} of $(M_t)_{t\in (0, T)}$ if there are sequences $\lambda_j\searrow 0,\,t_j\nearrow T$ and $(x_j)$ in $\rn$ such that
\begin{equation*}
(\Omega^j_s)_{s\in (a_j, b_j)}\to(\Omega'_s)_{s\in (-\infty, b)}
\end{equation*}
smoothly in compact subsets in space-time (that is in particular, the hypersurfaces $M^j_s=\partial\Omega^j_s$ converge smoothly). }\end{defi}

\begin{rem}{\rm \noindent For a solution
$(M_t)_{t\in [0, T)}$ which becomes singular for $t\nearrow T$, that is $\sup_{t<T}\sup_{M_t}|A|^2=\infty$ for the second fundamental form $A$ on $M_t$ one can
always find a rescaling limit for a suitable choice of sequences $(x_j)$ in $\rn$ and $(\lambda_j)\searrow 0$ (for example the reciprocal of the maximum of $|A|$
at an appropriately chosen sequence of times $t_j\nearrow T$). The smooth convergence follows from standard a priori estimates for mean curvature flow (see for
instance \cite{Hu2}).

Rescaling limits are so-called \emph {ancient solutions} which means that they have existed forever. Examples of ancient solutions are all homothetically
shrinking solutions of mean curvature flow such as the shrinking spheres given by $M'_s=\partial B_{\sqrt{-2ns}}$ for $s\in (-\infty, 0)$.

If the solution $(M_t)_{t\in (0, T)}$ has a so-called \emph{type II - singularity}, that is
\begin{equation*}
\sup_{t<T} \left((T-t)\sup_{M_t}|A|^2\right)=\infty ,
\end{equation*}
then by a rescaling process described in \cite{HS} one can even find a limit flow which is an \emph{eternal} solution, that is $b=\infty$. Examples of eternal
solutions are all \emph{stationary} solutions, that is solutions with $\Omega'_s =\Omega$ for all $s\in \mathbb{R}$. In this case, the hypersurface
$M=\partial\Omega$ is minimal that is satisfies $H =0$. Other eternal solutions are  \emph{translating} solutions of mean curvature flow for which
$\Omega'_s=\Omega+s \omega$ for $s\in\mathbb{R}$ where $\Omega\subset\rn$ and $\omega$ is a fixed unit vector in $\rn$. The corresponding hypersurfaces
$M=\partial\Omega$ satisfy the equation $H+\nu\cdot\omega = 0$. }\end{rem}

\bigskip

The statement of Corollary \ref{cor:mcf} is scaling invariant. Hence the rescaled solution $(M^j_s)_{s\in (a_j, b_j)}$ satisfies
\begin{equation*}
\frac{V(\Omega^j_s\cap B_r(x_0))}{r^{n+1}}\ge\kappa >0
\end{equation*}
for all $s\in (a_j, b_j)$ and $r\in (0, \sqrt{T}/\lambda_j)$  in balls with $V(\Omega^j_s\cap B_{r/2}(x_0))>0$ and
\begin{equation*}
\frac{V(\Omega^j_s\cap B_r(x_0))+r^2\int_{M^j_s\cap B_r(x_0)}|H|\,dS}{V(\Omega^j_s\cap B_{r/2}(x_0))}\le c_1.
\end{equation*}
The constant $\kappa = \kappa (n, \Omega_0, T, \sup_{M_0}|H|, c_1)$ is the same as for the unscaled solution. As a consequence we obtain a lower volume ratio
bound for rescaling limits, but without the radius restriction:

\begin{cor} Let $(M_t)_{t\in [0, T)}$ be a solution of mean curvature flow consisting of compact smooth, embedded hypersurfaces which enclose bounded regions
$(\Omega_t)_{t\in [0, T)}$ in $\rn$. Suppose furthermore that
inequality (C) holds. Then there is a constant $\kappa
>0$ depending only on $n, \Omega_0, T, \sup_{M_0}|H|$ and $c_1$
such that any rescaling limit $(M'_s)_{s\in (-\infty, b)}$ of
$(M_t)_{t\in [0, T)}$ with limiting enclosed regions
$(\Omega'_s)_{s\in (-\infty, b)}$ satisfies
\begin{equation*}
\frac{V(\Omega'_s\cap B_r(x_0))}{r^{n+1}}\ge\kappa
\end{equation*}
for every $s\in (-\infty, b)$ and $r >0$ in balls $B_r(x_0)$ with $V(\Omega'_s\cap B_{r/2}(x_0))>0$ and
\begin{equation*}
\frac{V(\Omega'_s\cap B_r(x_0))+r^2\int_{M'_s\cap B_r(x_0)}|H|\,dS}{V(\Omega'_s\cap B_{r/2}(x_0))}\le c_1.
\end{equation*}
\end{cor}

This Corollary rules out certain solutions of mean curvature flow
as rescaling limits under the assumption that our conjecture is
valid:

\begin{cor}\label{cor:mcflimits} If the conditions of the above corollary are satisfied then the following eternal solutions of mean curvature flow cannot occur
as rescaling limits of a compact, smooth embedded mean-convex
solution $(M_t)_{t\in [0, T)}$ of mean curvature flow which
encloses bounded regions $(\Omega_t)_{t\in [0, T)}$ in $\rn$:

\noindent (1) The stationary solution corresponding to a pair of
parallel hyperplanes that is given by $\Omega'_s =\Omega$ for all
$s\in\mathbb{R}$ where
\begin{equation*}
\Omega = \{x\in\rn,\,-d<x_{n+1}<d\}
\end{equation*}
for some $d>0$.

\noindent (2) The stationary solution of mean curvature flow
corresponding to the catenoid minimal surface $M=\partial\Omega$
in $\mathbb{R}^3$ given by
\begin{equation*}
\Omega =\{x=(\hat x, x_3)\in\mathbb{R}^3, \,|\hat x|\ge 1,\,|x_3|\le \cosh^{-1} |\hat x|\}.
\end{equation*}

\noindent (3) The translating solution corresponding to the grim reaper hypersurface $M=\partial\Omega$ where $\Omega=\mathbb{R}^{n-1} \times G$ with
\begin{equation*}
G=\left\{(x_n, x_{n+1})\in\mathbb{R}^2,\, -\pi/2 < x_n < \pi/2, \,x_{n+1}>-\log\cos  x_n\right\}.
\end{equation*}
\end{cor}

\noindent{\bf Proof.} All three examples admit sequences of balls
for radii increasing to infinity for which the volume ratio tends
to zero while the other quantities are controlled. This was
discussed in Corollary \ref{cor:-infty}.

\bigskip
\begin{rem}{\rm \noindent (1) In the special situation where the original solution $(M_t)$ is mean convex, that is $H>0$ for $M_0$ and subsequently for all
$M_t$ by the maximum principle,
White [Wh] ruled out the grim reaper hypersurface as a rescaling limit using techniques from minimal surface theory and geometric measure theory. His methods
extend also to non-smooth limit flows of generalized mean curvature flow solutions in the mean-convex case.

(2) In view of Corollary \ref{cor:-infty}, the first two examples satisfy $\inf_{\tau >0}\mu(\Omega,\tau)=-\infty$ and the third one $\inf_{\tau
>0}\mu_H(\Omega,\tau)=-\infty$.

(3) The embeddedness assumption on the hypersurfaces $M_t$ is essential. In \cite{A}, it is proved that rescaling limits of non-embedded planar curves near
singularities are given by the grim reaper curve $\Gamma =\partial G$ defined above.

(4) Some other translating solutions can occur as rescaling limits such as for instance a rotationally symmetric translating bowl (see for instance \cite{Wa}).
The region bounded by this translating bowl opens up quadratically so one can show that it satisfies the conclusions of the above corollary.

(5) For the shrinking solution $\Omega'_s=B_{\sqrt{-2ns}}$ there is no lower bound of the form
\begin{equation*}
\frac{V(\Omega'_s\cap B_r)}{r^{n+1}}\ge\kappa >0
\end{equation*}
with a fixed $\kappa$ for all $s<0$ and all $r>0$ since the balls $\Omega'_s$ shrink to the origin for $s\nearrow 0$. This does not contradict the corollary
though as $\kappa$ depends on $c_1$ and in this case $c_1$ behaves like $-c(n)r^2s^{-1}$ since for $s\in [-1/(2n), 0)$ and $r\ge 1$
\begin{equation*}
\frac{r^2\int_{M'_s\cap B_r}H\,dS}{V(\Omega'_s\cap B_{r/2})}=\frac{\int_{M'_s}H\,dS}{V(\Omega'_s)}=-c(n)\frac{r^2}{s}.
\end{equation*}}
\end{rem}
\bigskip

\section*{Appendix A. Some basic properties of entropies in $\rn$}\label{examples}

In this appendix, we discuss some explicit examples of entropies in $\rn$.

\bigskip
\noindent  (1) When $\beta =0$ and $\Omega =\rn$ we have (see \cite{P})
\begin{equation*}
\mathcal{W}(\rn, f, \tau)=\int_{\rn}\left(\tau|\nabla f|^2+f - (n+1)\right) u\, dx\ge 0
\end{equation*}
for all $f$ satisfying
\begin{equation*}
\int_{\rn}u \,dx =1
\end{equation*}
with equality when $f(x)=\frac{|x|^2}{4\tau}$. In particular therefore
\begin{equation*}
\mu (\rn, \tau) =0
\end{equation*}
for all $\tau >0$.

This is the Gaussian logarithmic Sobolev inequality due to L.Gross (\cite{G}). Scaling by $x=\sqrt{2\tau}y$, setting $f=\frac{|y|^2}{2}-\log\varphi^2$ as in
\cite{P} and using the identity
\begin{equation*}
\int_{\rn}\left(|y|^2-(n+1)\right)\gamma_{n+1}\,dy=-\int_{\rn}{\rm div}\left(y\gamma_{n+1}\right)\,dy=0
\end{equation*}
for the Gaussian
\begin{equation*}
\gamma_{n+1}(y)=\frac{e^{-\frac{|y|^2}{2}}}{({2\pi})^{n+1}}
\end{equation*}
we obtain its standard form
\begin{equation*}
\int_{\rn}\varphi^2\log\varphi\;\gamma_{n+1}\,dy\le\frac{1}{2}\int_{\rn}|\nabla\varphi|^2\;\gamma_{n+1}\,dy
\end{equation*}
for all $\varphi$ satisfying
\begin{equation*}
\int_{\rn}\varphi^2\gamma_{n+1}\,dx=1.
\end{equation*}

\noindent (2) For $x\in\Omega\subset\rn$ we set $x=\lambda y+x_0$ where $\lambda >0$ and $x_0\in\rn$. We then obtain
\begin{align}
\mathcal{W}_\beta (\Omega, f, \tau)&=
\mathcal{W}(\lambda^{-1}(\Omega -x_0), f(\lambda\,\cdot+\;x_0), \lambda^{-2}\tau)\nonumber\\
&\qquad\qquad+ 2(\lambda^{-2}\tau)\int_{\lambda^{-1}(\partial\Omega -x_0)}\lambda\beta (\lambda y+x_0)\frac{e^{-f(\lambda
y+x_0)}}{\left(4\pi\lambda^{-2}\tau\right)^{\frac{n+1}{2}}}\,dS (y)\nonumber
\end{align}
and
\begin{equation*}
1=\int_\Omega u(x)\,dx=\int_{\lambda^{-1}(\Omega -x_0)} u(\lambda y+x_0)\,dy.
\end{equation*}
Therefore
\begin{equation*}
\mu_\beta(\Omega, \tau)=\mu_{\lambda\beta (\lambda\,\cdot+x_0)}(\lambda^{-1}(\Omega -x_0),\lambda^{-2}\tau).
\end{equation*}

\noindent Suppose that $\beta :\rn\to \mathbb{R}$ satisfies
\begin{equation*}
\beta\left (\frac {x-x_0}{\lambda}\right)=\lambda\beta(x)
\end{equation*}
for $x, x_0\in\rn$ and $\lambda >0$ or that $\beta =\beta_{\partial\Omega}$ is a geometric quantity which behaves like
\begin{equation*}
\beta (y)=\lambda\beta (x)
\end{equation*}
where $x=\lambda y +x_0\in\partial\Omega$ for $y\in\frac{1}{\lambda}(\partial\Omega -x_0)$ such as for example the mean curvature of $\partial\Omega$. Then
\begin{equation*}
\mu_\beta(\Omega, \tau)=\mu_\beta(\lambda^{-1}(\Omega -x_0),\lambda^{-2}\tau).
\end{equation*}

\noindent For $x_0=0, \lambda =\sqrt{2\tau}, \Omega$ replaced by $\sqrt{2\tau}\,\Omega$ and such functions $\beta$ this yields
\begin{equation*}
\mu_\beta(\sqrt{2\tau}\,\Omega, \tau)=\mu_\beta(\Omega,1/2).
\end{equation*}

\noindent (3) If $x_0\in\Omega$ then
\begin{equation*}
\lambda^{-1}(\Omega -x_0)\to\rn.
\end{equation*}
Using this, the scaling identity for $\mu_\beta$ with $\lambda = \sqrt{2\tau}$ as well as the identity $\mu(\rn, 1/2)=0$ we expect that
\begin{equation*}
\mu_\beta(\Omega, \tau)\to 0
\end{equation*}
for $\tau\to 0$. This should follow along the same lines as in \cite{N}.

\noindent (4) A natural example is
\begin{equation*}
\beta=\frac{x\cdot\nu}{2\tau}
\end{equation*}
where $\nu$ is the unit outward pointing normal to $\partial\Omega$. By the above scaling property we have
\begin{equation*}
\mu_{\frac{x\cdot\nu}{2\tau}}(\sqrt{2\tau}\,\Omega, \tau)  =\mu_{y\cdot\nu}(\Omega,1/2)
\end{equation*}
where $x=\sqrt{2\tau} y$ and $y\in\Omega$.

An example of a function $f$ on $\Omega\subset\rn$ satisfying the normalisation condition
\begin{equation*}
\int_\Omega\frac{e^{-f}}{(4\pi\tau)^{\frac{n+1}{2}}}\,dx=1
\end{equation*}
is
\begin{equation*}
f=\frac{|x|^2}{4\tau}-\log c
\end{equation*}
where
\begin{equation*}
\frac{1}{c}=\int_\Omega \frac{e^{-\frac{|x|^2}{4\tau}}}{(4\pi\tau)^{\frac{n+1}{2}}}\,dx.
\end{equation*}
For this $f$ and $\beta =\frac{x\cdot\nu}{2\tau}$ one calculates
\begin{equation*}
\mathcal{W}_\beta (\Omega, f, \tau)=c\left(\int_\Omega\left(\frac{|x|^2}{2\tau}-
(n+1)\right)\frac{e^{-\frac{|x|^2}{4\tau}}}{(4\pi\tau)^{\frac{n+1}{2}}}\,dx+\int_{\partial\Omega}x\cdot\nu
\frac{e^{-\frac{|x|^2}{4\tau}}}{(4\pi\tau)^{\frac{n+1}{2}}}\,dS\right) +\log c .
\end{equation*}
Since
\begin{equation*}
{\rm div}\left( xe^{-\frac{|x|^2}{4\tau}}\right)=-\left(\frac{|x|^2}{2\tau}- (n+1)\right) e^{-\frac{|x|^2}{4\tau}}
\end{equation*}
this implies
\begin{equation*}
\mathcal{W}_\beta (\Omega, f, \tau)=\log c
\end{equation*}
by the divergence theorem.

Note that for $\Omega = \rn$ we have $c=1$ and hence $\mathcal{W}_\beta (\Omega, f, \tau)=\mathcal{W}(\Omega, f, \tau)=0$ for $f=\frac{|x|^2}{4\tau}$.

For the half-space $H_a=\{x\in\rn, x_{n+1}<a\}, \,a\in\mathbb{R}$ and $\beta =\frac{x\cdot\nu}{2\tau}$ we calculate
\begin{equation*}
\frac{1}{c}=\int_{-\infty}^{\frac{a}{\sqrt{2\tau}}}e^{-\frac{z^2}{2}}\,dz.
\end{equation*}
This implies that $\mu_{\frac{x\cdot\nu}{2\tau}}(H_a, \tau)\to -\infty$ for $a\to -\infty$ as well as $\lim_{\tau\to 0}\mathcal{W}_{\frac{x\cdot\nu}{2\tau}}
(H_a, f, \tau)=0$ and $\lim_{\tau\to 0}\mathcal{W}_{\frac{x\cdot\nu}{2\tau}} (H_a, f, \tau)=-\log 2 <0$ for fixed $a\in\mathbb{R}$.

By the scaling and translation property above we have
\begin{equation*}
\mathcal{W}_{\frac{x\cdot\nu}{2\tau}}\left(\Omega, \frac{|x|^2}{4\tau}-\log c, \tau\right)=\mathcal{W}_{y\cdot\nu} \left(\frac{1}{\sqrt{2\tau}}\Omega,
\frac{|y|^2}{2}-\log c, \frac{1}{2}\right)=\log c
\end{equation*}
for $x=\sqrt{2\tau}y\in\Omega$ with the condition
\begin{equation*}
\int_{\frac{1}{\sqrt2\tau}\Omega}\gamma_{n+1}\,dy =\frac{1}{c}.
\end{equation*}
If the $(n+1)$ -dimensional volume of a set $\Omega$ inside large balls grows like
\begin{equation*}
V(\Omega\cap B_R)\le c R^p
\end{equation*}
for $R\ge R_0$ and $p<n+1$ one checks that
\begin{equation*}
\int_{\frac{1}{\sqrt2\tau}\Omega}\gamma_{n+1}\,dy\to 0
\end{equation*}
for $\tau\to\infty$. Therefore $c\to\infty$ and hence
\begin{equation*}
\mu_{\frac{x\cdot\nu}{2\tau}}(\Omega, \tau)\to -\infty.
\end{equation*}
Such sets $\Omega$ include for instance all bounded sets but also unbounded sets which lie in a slab in $\rn$. In the latter case the volume in balls grows like
$R^n$.

\bigskip

\section*{Appendix B. Sobolev and logarithmic Sobolev inequalities}\label{sobolev}

For the convenience of the reader who is unfamiliar with logarithmic Sobolev inequalities we show how these can be derived from the standard Sobolev inequality.
We essentially follow a proof given in \cite{G}.

\begin{unnumberedtheorem}[Logarithmic Sobolev inequality]
For any open subset $\Omega$ of a Riemannian manifold $(X, g)$ which satisfies the Sobolev inequality
\begin{equation*}
\left(\int_\Omega |\psi|^{\frac{n+1}{n}}\,dV\right)^{\frac{n}{n+1}}\le c_S(\Omega, g)\int_\Omega\left(|\nabla\psi|+|\psi|\right)\,dV
\end{equation*}
for all $\psi\in C^1(\bar\Omega)$ there also holds a logarithmic Sobolev inequality of the form
\begin{equation*}
\int_\Omega\left(\epsilon\; |\nabla\varphi|^2 -\varphi^2\log\varphi^2\right)\,dV\ge -c(n)(1+\log c_S(\Omega, g)) - \frac{1}{\epsilon}
\end{equation*}
for functions $\varphi$ satisfying $\int_\Omega \varphi^2\,dV =1$ and every $\epsilon >0$.
\end{unnumberedtheorem}

\noindent{\bf Proof.} By a standard approximation argument it will be sufficient to prove the theorem for non-negative functions. We abbreviate
\begin{equation*}
\Vert\psi\Vert_p \equiv \left(\int_\Omega \psi^p\,dV\right)^{\frac{1}{p}}
\end{equation*}
for $p>0$. The interpolation inequality for $L^p$ -norms says for functions $\psi$ satisfying $\Vert\psi\Vert_1=1$ that
\begin{equation*}
\Vert\psi\Vert_q\le \Vert\psi\Vert_{\frac{n}{n-1}}^{n-\frac{n}{q}}
\end{equation*}
for $1\le q\le \frac{n}{n-1}$. Since for $q=1$ we have equality, differentiation with respect to $q$ at $q=1$ preserves the inequality and leads to
\begin{equation*}
\int_\Omega\psi\log\psi\,dV\le n\log\Vert\psi\Vert_{\frac{n}{n-1}}.
\end{equation*}
In view of the Sobolev inequality
\begin{equation*}
\Vert\psi\Vert_{\frac{n}{n-1}}\le c_S(\Omega)\left(\Vert\nabla\psi\Vert_1+1\right)
\end{equation*}
for such functions this yields
\begin{align}
\int_\Omega\psi\log\psi\,dV&\le n\log\left(c_S(\Omega)\left(\Vert\nabla\psi\Vert_1+1\right)\right)\nonumber\\
&=n\log\left(\frac{1}{n}\left(\Vert\nabla\psi\Vert_1+1\right)\right)+n\log
(n c_S(\Omega)). \nonumber
\end{align}
The inequality $\log x\le x-1$ implies
\begin{equation*}
\int_\Omega\psi\log\psi\,dV\le\Vert\nabla\psi\Vert_1+c(n)(1+\log c_S(\Omega)).
\end{equation*}
Setting $\psi=\varphi^2$ with $\int_\Omega\varphi^2\,dV=1$ gives
\begin{equation*}
\int_\Omega\varphi^2\log\varphi^2\,dV\le\int_\Omega |\nabla\varphi^2|\,dV+c(n)(1+\log c_S(\Omega)).
\end{equation*}
Using Young's inequality, we finally arrive
\begin{equation*}
\int_\Omega\varphi^2\log\varphi^2\,dV\le\epsilon\int_\Omega |\nabla\varphi|^2\,dV+\frac{1}{\epsilon}+c(n)(1+\log c_S(\Omega))
\end{equation*}
where we again used $\int_\Omega\varphi^2\,dV=1$. \qed

\bigskip
\section*{Appendix C. Proof of the Evolution equation for $W$}\label{w}

For the convenience of the reader we give a detailed proof of the
evolution equation of Proposition \ref{pro:wequation} in Section
\ref{ricci}. In Section \ref{ricci}, we merely modified the
appropriate formulas in \cite{P} and \cite{N} by transforming to
total time derivatives.

Let us briefly recall the set-up given in Section \ref{ricci} in the case of evolving domains in $\rn$. We consider a family of subsets $(\Omega_t)_{t\in (0,
T)}$ in $\rn$  which evolve by the equation
\begin{equation}\label{eq:omegaflowa}
\frac{\partial x}{\partial t}=-\nabla f(x, t)
\end{equation}
for $x\in\Omega_t$ where $f(t)$ satisfies the equation
\begin{equation}\label{eq:fequationa}
\left(\frac{\partial }{\partial t} +\Delta\right) f= |\nabla f|^2+\frac{n+1}{2\tau}
\end{equation}
in $\Omega_t$ for $t\in (0, T)$. The total time derivative  of $f$ is given by
\begin{equation}\label{eq:totaldera}
\frac{df}{dt}=\frac{\partial f}{\partial t}+\left\langle\nabla f, \frac{\partial x}{\partial t}\right\rangle=\frac{\partial f}{\partial t}-|\nabla f|^2
\end{equation}
and so (\ref{eq:fequationa}) can also be written as
\begin{equation}\label{eq:totalfequationa}
\left(\frac{d}{dt} +\Delta\right) f= \frac{n+1}{2\tau}.
\end{equation}
We also assume that $\tau(t)>0$ evolves by $\frac{\partial\tau}{\partial t}=-1$.

\bigskip

\noindent {\bf Proposition.} \emph{In the above setting, the function $W=\tau(2\Delta f-|\nabla f|^2)+f-(n+1)$ satisfies the evolution equation
\begin{equation*}
\left(\frac{d}{dt} +\Delta\right)W= 2\tau \left|\nabla_i\nabla_j f-\frac{1}{2\tau}\,\delta_{ij}\right|^2 + \nabla W\cdot\nabla f.
\end{equation*}}

\noindent{\bf Proof.} We adapt the computation in \cite{N} to the case of domains evolving by (\ref{eq:omegaflowa}) (the different sign in Ni's Lemma 2.2 stems
from the fact that he considers the forward heat equation by interchanging the roles of $\tau$ and $t$.) In a general Riemannian manifold $(X, g)$ an additional
Ricci term arises when we interchange third derivatives of $f$. In the Ricci flow case this expression is balanced by terms coming from the time derivative of
the metric. Details of the latter can be found in \cite{KL}.

If we write above $x=\phi(q, t)$ where $\phi_t =\phi (\cdot,
t):\Omega\to\Omega_t$ are the diffeomorphisms evolving $\Omega_t$,
 the pulled back function $f$ given by
\begin{equation*}
\tilde f(q, t)=f (\phi(q, t), t)
\end{equation*}
satisfies
\begin{equation*}
\frac{df}{dt}(x, t)=\frac{\partial\tilde f}{\partial t}(q, t).
\end{equation*}
The evolution equation (\ref{eq:omegaflowa}) written in terms of
$\tilde f(q, t)=f(\phi(q, t), t)$ looks like
\begin{equation*}
\frac{\partial\phi}{\partial t}(q, t)=-\tilde\nabla\tilde f (q ,t)
\end{equation*}
where $\tilde\nabla$ is the gradient with respect to the pull-back
of the Euclidean metric under $\phi_t$ on $\Omega\subset\rn$ given
by
\begin{equation*}
g_{ij}(q, t)=\frac{\partial\phi}{\partial q_i}(q,
t)\cdot\frac{\partial\phi}{\partial q_j}(q, t).
\end{equation*}
In these coordinates we have
\begin{equation*}
|\nabla f|^2=g^{ij}\frac{\partial\tilde f}{\partial q_i}\frac{\partial\tilde f}{\partial q_j}.
\end{equation*}
One now calculates
\begin{equation*}
\frac{\partial}{\partial t}g_{ij}=-2\tilde\nabla_i\tilde\nabla_jf,
\end{equation*}
and the inverse metric satisfies
\begin{equation*}
\frac{\partial}{\partial t}g^{ij}=2\tilde\nabla^i\tilde\nabla^jf.
\end{equation*}
Furthermore one computes for the Christoffel symbols of the $g_{ij}$
\begin{equation}\label{eq:christoffel}
g^{ij}\frac{\partial}{\partial t}\Gamma_{ij}^k=\tilde\nabla^k\tilde\Delta \tilde f.
\end{equation}
One then checks from this and $\Delta f(x, t)=\tilde\Delta\tilde f(q, t)$ with
\begin{equation*}
\tilde\Delta\tilde f= g^{ij}\left(\frac{\partial^2\tilde f}{\partial q_i\partial q_j}-\Gamma_{ij}^k\frac{\partial\tilde f}{\partial q_k}\right)
\end{equation*}
that the identities
\begin{equation}\label{eq:dtgrad}
\frac{d}{dt}|\nabla f|^2=2\nabla_i\nabla_jf\nabla_if\nabla_jf+2\nabla f\cdot\nabla\frac{df}{dt}
\end{equation}
and
\begin{equation}\label{eq:dtdeltaf}
\frac{d}{dt}\Delta f=\Delta\frac{df}{dt}+2|\nabla^2f|^2+\nabla f\cdot\nabla\Delta f
\end{equation}
hold. We now follow \cite{N} exactly, except for working with $\frac{df}{dt}$ instead of $\frac{\partial f}{\partial t}-|\nabla f|^2$. The latter of the above
identities in combination with (\ref{eq:totalfequationa}) and the relation $\frac{\partial\tau}{\partial t}=-1$ implies
\begin{equation}\label{eq:dfdtequation}
\left(\frac{d}{dt} +\Delta\right)\frac{df}{dt}=\frac{n+1}{2\tau^2}-2|\nabla^2f|^2-\nabla f\cdot\nabla\Delta f.
\end{equation}
Combining (\ref{eq:dtgrad}) and (\ref{eq:totalfequationa}) with the Bochner identity
\begin{equation*}\Delta |\nabla f|^2=2|\nabla^2f|^2+2\nabla f\cdot\nabla\Delta f
\end{equation*}
we find
\begin{equation}\label{eq:gradevolutionequation}
\left(\frac{d}{dt} +\Delta\right)|\nabla f|^2=2|\nabla^2f|^2 +2\nabla_i\nabla_jf\nabla_if\nabla_jf.
\end{equation}

To break up the calculation for $W$, we rewrite $W=\tau (2\Delta f -|\nabla f|^2)+f-(n+1)$ using (\ref{eq:totalfequationa}) as
\begin{equation*}
W=\tau w+f
\end{equation*}
where
\begin{equation}\label{eq:wdefinition}
w=-2\frac{df}{dt}-|\nabla f|^2=2\Delta f-|\nabla f|^2-\frac{n+1}{\tau}.
\end{equation}
From (\ref{eq:dfdtequation}) and (\ref{eq:gradevolutionequation}) we calculate
\begin{equation*}
\left(\frac{d}{dt} +\Delta\right)w=2|\nabla^2f|^2-\frac{n+1}{\tau^2}-2\nabla_i\nabla_jf\nabla_if\nabla_jf+2\nabla f\cdot\nabla\Delta f.
\end{equation*}
Since
\begin{equation*}
-2\nabla_i\nabla_jf\nabla_if\nabla_jf+2\nabla f\cdot\nabla\Delta f=\nabla f\cdot\nabla w
\end{equation*}
we thus arrive at
\begin{equation*}
\left(\frac{d}{dt} +\Delta\right)w=2|\nabla^2f|^2-\frac{n+1}{\tau^2}+\nabla f\cdot\nabla w.
\end{equation*}
Using again $\frac{\partial\tau}{\partial t}=-1$ we now compute
\begin{align}
\left(\frac{d}{dt} +\Delta\right)W&=\left(\frac{d}{dt} +\Delta\right)(\tau w+f)\nonumber\\
&=-w+2\tau |\nabla^2f|^2-\frac{n+1}{\tau}+\nabla f\cdot\nabla (\tau w) +\frac{n+1}{2\tau}\nonumber\\
&=\nabla f\cdot\nabla W-w-|\nabla f|^2+2\tau |\nabla^2f|^2-\frac{n+1}{2\tau}.\nonumber
\end{align}
Substituting the identities
\begin{equation*}
2\tau |\nabla^2f|^2=2\tau\left|\nabla_i\nabla_jf-\frac{\delta_{ij}}{2\tau}\right|^2+2\Delta f -\frac{n+1}{2\tau}
\end{equation*}
and
\begin{equation*}
w=2\Delta f-|\nabla f|^2-\frac{n+1}{\tau}
\end{equation*}
yields the desired evolution equation for $W$.\qed

{\small\small \noindent Fachbereich Mathematik und Informatik\\
Freie Universit\"at Berlin\\ Arnimallee 2-6\\ 14195
Berlin\\Germany\\}

\end{document}